\documentclass[10pt]{article}
\usepackage{amsmath}
\usepackage{amsfonts}
\usepackage{amssymb}
\usepackage{comment}
\newtheorem{theorem}{Theorem}[section]
\newtheorem{definition}[theorem]{Definition}
\newcommand{\Z}{\mathbb Z}
\newcommand{\Q}{\mathbb Q}
\newcommand{\G}{\Gamma}

\newcommand{\be}{\beta}

\newcommand{\z}{\zeta}

\begin{document}
\pagestyle{plain}
\title{A Database of Continued Fractions of Polynomial Type}
\author{Henri Cohen}

\maketitle

\begin{abstract}
  We describe a database of 1925 continued fractions with polynomial
  coefficients, of which more than 1600 are new, both for interesting
  constants and for transcendental functions, and provide the database inside
  the \TeX\ source of the paper. Look in particular at the section
  ``Table of Contents''.
\end{abstract}  
  
\section{Introduction}

{\bf Please download the \TeX\ source of this paper}
since it contains the database described here, otherwise this paper
is useless.

\smallskip

This is not quite a standard mathematical paper, but I believe it will be
useful (and hopefully become a standard reference) to any mathematician
interested in continued fractions, in particular attached to interesting
constants or functions.

I have developed an extensive {\tt Pari/GP} package for dealing with
continued fractions (abbreviated CF in the sequel), for now accessible only
as a GIT branch under the name {\tt henri-ellCF3} (please read the
{\tt Pari/GP} installation documentation to obtain this). Thanks to
this package, I have done several things.

\begin{enumerate}\item First, starting from the 250 or so CFs with
  polynomial coefficients (see below for the precise definition) existing in
  the literature, I have succeeded in creating more than 1650 additional ones,
  some of them I believe quite interesting.
\item Second, for almost \emph{every} CF in the corresponding database, I
  have computed the precise rate of convergence of the partial quotients.
  This is almost never done in the literature, although many texts deal
  with convergence of CFs. To give an almost random example, we have
  $$\log(2)=3/(5-9\cdot1^2/(15-9\cdot2^2/(25-9\cdot3^2/(35-9\cdot4^2-\cdots))))\;,$$
  and in addition to this CF the database contains the following:
  $$\log(2)-\dfrac{p(n)}{q(n)}=\dfrac{\pi/3}{9^n}\left(1-\dfrac{5/16}{n}+\dfrac{105/512}{n^2}-\dfrac{1725/8192}{n^3}+\cdots\right)\;.$$
\item Third, in many cases, I have given some explanation on the origin
  of the CF, in particular when it has been obtained by Ap\'ery acceleration
  techniques, see \cite{Coh2}.
\item Fourth, in many cases the CF can also be expressed as a hypergeometric
  series, which is given.
\end{enumerate}

\medskip

I have created a file containing these 1925 continued fractions
together with their relevant information, which is readable from {\tt Pari/GP}
with the command {\tt V=readvec("file.gp")}, if {\tt file.gp} is the name
given to the file (it is of course also readable by other programs).
To obtain this file, from the arXiv repository from which
you downloaded the pdf you must also obtain the \TeX\ source, and the file
is entirely contained between the last two {\tt comment} commands.

\medskip

Special thanks go to Francesco Del Re who used AI to go though both the printed dictionary and the
present GP database and corrected a large number of mistakes and typos.

\section{The Database in Detail}

I will now describe the precise format of the database of continued fractions.
It is intended to be used by {\tt Pari/GP}, but of course since the format
is transparent it can be used with any other software through an appropriate
interface.

\smallskip

\subsection{Computer Representation of Continued Fractions}

\smallskip

Let $(a(n))_{n\ge0}$ and $(b(n))_{n\ge0}$ be two sequences of complex numbers.
The continued fraction $(a,b)$ is the formal expression
$$(a,b)=a(0)+\dfrac{b(0)}{a(1)+\dfrac{b(1)}{a(2)+\dfrac{b(2)}{a(3)+\ddots}}}$$
written more compactly as $(a,b)=a(0)+b(0)/(a(1)+b(1)/(a(2)+\cdots))$.
Note that this notation is different from that used by other authors
such as \cite{Cuyt}. I do not claim that my notation
is better but simply that it is the one used in the present paper.
We will always assume that $b(n)\ne0$ for all $n$ (otherwise the CF
terminates), and that $a(n)=0$ only for a finite number of $n$. We write as
usual formally $p(n)/q(n)=a(0)+b(0)/(a(1)+b(1)/(a(2)+\cdots+b(n-1)/a(n)))$.

\begin{definition} A CF $((a(n))_{n\ge0},(b(n))_{n\ge0})$ is said to have
  \emph{polynomial coefficients} if there exists a \emph{period} $T$ such
  that for every $j$ with $0\le j<T$, for $n$ sufficiently large both
  $a(nT+j)$ and $b(nT+j)$ are rational functions of $n$ with
  \emph{rational} coefficients.\end{definition}

The name \emph{polynomial} is justified since it is clear that
any CF with coefficients in $\Q(n)$ is equivalent (i.e., a CF with the
same $p(n)/q(n)$) to one with coefficients in $\Z[n]$.

In practice $T$ will be equal to $1$ or $2$, so for $n$ sufficiently large
either $a(n)$ and $b(n)$ are rational functions, or $a(2n)$, $a(2n+1)$,
$b(2n)$, and $b(2n+1)$ are rational functions.

It is evidently easy to \emph{contract} a CF of period 2 to a CF of period
1, but the point of keeping a period 2 CF is that often its coefficients
are simple, while those of the contracted CF are complicated.

\smallskip

Note that (of course excluding finite and periodic CFs) most explicit CFs
occurring in the literature have polynomial coefficients,
with the notable exception of $q$-continued fractions such as Ramanujan's
$$R(q)=q^{1/5}/(1+q/(1+q^2/(1+q^3/(1+\cdots))))\;,$$
which are fascinating but will not be considered here.

\medskip

Let $(a,b)=((a(n))_{n\ge0},(b(n))_{n\ge0})$ be a continued fraction with
polynomial coefficients. If $a$ (similarly $b$) has period $1$, i.e., if
$a(n)$ is a polynomial for $n\ge n_0$ for some $n_0$, we will represent $a$
as the vector $$A=[a(0),a(1),\dotsc,a(n_0-1),a(n)]\;,$$ where for
$i\le n_0-1$ the $a(i)$ do not depend on $n$, while $a(n)$ is a polynomial
(exceptionally also a rational function) in the reserved variable $n$. Note
that the variable $n$ is \emph{reserved} for this purpose and cannot be used
elsewhere. It would perhaps have been more logical to use the variable $x$
instead, but the variable $n$ is almost universally used in texts so I have
chosen to use it.

Note also that the representation is not unique, for instance
we also have $A=[a(0),a(1),\dotsc,a(n_0-1),a(n_0),a(n)]$.

For instance, the continued fraction
$$\log(2)=3/(5-9\cdot1^2/(15-9\cdot2^2/(25-9\cdot3^2/(35-9\cdot4^2-\cdots))))$$
given above is represented by

\centerline{\tt [[0,10*n-5],[3,-9*n\^{}2]]}

\smallskip

If $A$ has period $2$, i.e., if $a(2n)$ and $a(2n+1)$ are possibly different
polynomials for $n$ large, we will represent $A$ as a vector of $2$-component
vectors $$A=[[a(0),a(1)],\dotsc,[a(2m_0-2),a(2m_0-1)],[a(2n),a(2n+1)]]\;.$$

For instance, for the continued fraction
$$\pi=4-4/(3+9/(5+4/(7+25/(9+16/(11+49/(13+36/(15+\cdots)))))))\;,$$
$A$ has period 1 but $B$ has period 2, so we can represent the CF by

\centerline{\tt [[4,2*n+1],[[-4,9],[(2*n)\^{}2,(2*n+3)\^{}2]]]}

Note that in this case \emph{all} the components of $B$ must be two-component
vectors. For instance the above could \emph{not} be written as

\centerline{\tt [[4,2*n+1],[-4,[(2*n+1)\^{}2,(2*n)\^{}2]]]}

\medskip

\subsection{The Database Formats}

\medskip

First, a warning. We reserve the following variable names which must never
be used in working with the CFs of this file:

\begin{enumerate}\item ``n'', which is universally used to denote the CF index.
\item ``d'', which denotes the radicand of a quadratic extension, see below
  the explanation for {\tt v[3][1]}.
\item ``x'' which stands for $1/n$ and ``y'' which stands for $1/n^{1/2}$
  in the asymptotic expansions $A$, see the explanation for {\tt v[3][4]}.
\end{enumerate}

\smallskip

There are three kinds of inputs in the file. By far the most important ones are
of course the CFs themselves, which are also the most numerous.

\smallskip

$\bullet$ {\bf Format of the Continued Fractions Entries}

Such an entry, let's call it {\tt v}, has eight components, which we now
describe. To make the explanation clearer, we use the following example:

\begin{verbatim}
[(z)->cosh(Pi*z),[[1,1/2,4*n^2-8*n+(2*z^2+5)],
[2*z^2,-4*n^4+8*n^3+(-4*z^2-6)*n^2+(4*z^2+2)*n+(-z^2-1/4)]],
[1,[0,1,0,1],(z)->z^2*cosh(Pi*z),
1-1/2*z^2*x+(1/6*z^4-1/6*z^2-1/12)*x^2+(-1/24*z^4+1/6*z^2+1/12)*x^3+O(x^4)],
[[1,4*z^2,[1/2-I*z,1;1/2+I*z,1;3/2,-2],1],
 [1,-4*z^2,[-1/2,2;1/2-I*z,-1;1/2+I*z,-1],1]],
[],"6.6.5",["AP->6.6.11","APD->6.6.6"],""];
\end{verbatim}

\begin{enumerate}\item {\tt v[1]} is a closure giving the limit of the CF. Note
  that in CFs involving variables, the CF may converge to that limit only
  in a certain range of the variables, but we have not given that (we may
  in a future version).
\item {\tt v[2]} is the $[a(n),b(n)]$ of the CF in the format
  explained above, with one important difference: consider for instance the
  entry 5.7.1. beginning with:
\begin{verbatim}
  [k->zeta(k),[0,n^k+(n-1)^k],[1,-n^(2*k)]]
\end{verbatim}
Giving this to {\tt GP} would make it choke (give an error message) because
{\tt n\^{}k} is not a valid {\tt GP} input when {\tt k} is a variable.
  Thus we replace it in the file by
\begin{verbatim}
  [k->zeta(k),[k->[0,n^k+(n-1)^k],k->[1,-n^(2*k)]]]
\end{verbatim}
which is a valid input since the closure is not evaluated.
\item {\tt v[3]} is itself a 4-component vector {\tt [d,[F,E,D,P],C,A]}
  giving the exact speed of convergence of the CF, including a few terms of
  its asymptotic expansion. If $L$ denotes the value of the CF and $p(n)/q(n)$
  are the partial quotients, this notation means that
  $$L-\dfrac{p(n)}{q(n)}=\dfrac{C}{n!^FE^ne^{\sqrt{Dn}}n^P}A\;,$$
  with the following conventions:
  \begin{enumerate}
  \item {\tt v[3][1]} is a number or a polynomial (very often 1) with the
    following meaning: often, the speed of convergence and the asymptotic
    expansions involve elements of a quadratic extension (of $\Q$ or of
    $\Q[z]$ for instance). To avoid using square roots, we set {\tt v[3][1]}
    to the radicand, and we use the specific letter {\tt d} to indicate the
    square root of {\tt v[3][1]}. Thus, we must avoid the letter {\tt d} in
    variables of CFs. In the running example, we have {\tt v[3][1]=1} so
    no quadratic extension.
  \item {\tt v[3][2]} gives the speed of convergence in the above FEDP format,
    using {\tt d} if necessary. The running example has
    {\tt v[3][2]=[0,1,0,1]}, meaning that we are in case $P^+$ with
    $P=1$. In some rare cases, the speed does not make any sense, in which
    case we set {\tt v[3][2]=0}. In some other rare
    cases, the CF does not converge for real values of the variables (which
    is implicitly what is assumed), in which case we set
    {\tt v[3][2]=[0,0,0,0]}.
  \item {\tt v[3][3]} gives the constant $C$ such that
    $S-p(n)/q(n)\sim C/[FEDP]$, or $0$ if unknown or has not been computed.
    Note that $C$ is always given as a closure, first because in some cases
    (as for the CF itself) it is not possible to do otherwise in {\tt GP},
    and second because writing {\tt C=()->Pi} is both much clearer to read
    than {\tt C=3.1415...} and also allows to have {\tt C} to any accuracy.
  \item {\tt v[3][4]} gives the asympotic expansion {\tt A=1+...}, such that
    $S-p(n)/q(n)=(C/[FEDP])A$, or $0$ if unknown or not computed.
    The reserved letter {\tt x} means $1/n$ and the reserved letter {\tt y}
    means $1/n^{1/2}$.

    Thus, the above example says that
    $$\cosh(\pi z)-\dfrac{p(n)}{q(n)}=\dfrac{z^2\cosh(\pi z)}{n}\left(1-\dfrac{z^2/2}{n}+\dfrac{z^4/6-z^2/2-1/12}{n^2}+\cdots\right)$$
  \end{enumerate}
\item In a large number of cases, the CF or its inverse corresponds to a more
  or less explicit hypergeometric series {\tt S} given in a special format
  explained below. Thus, {\tt v[4]} is set to {\tt 0} if neither CF nor
  its inverse has such a series, to {\tt S} if only the CF has a series,
  to the two-component vector {\tt [S,T]} if both have a series (with {\tt T}
  the series for the inverse), and finally the two-component vector
  {\tt [[],T]} if only the inverse has a series.

  These series are given in the following format: they are all of the form
  $$S=u+\sum_{n\ge0}h(n)\prod_{1\le i\le g}(a_i)_n^{e_i}w^n\;,$$
  where $h$ is a rational function, $e_i\in\Z$, and as usual $(a_i)_n$ denotes
  the rising Pochhammer symbol $\G(a_i+n)/\G(a_i)$ (note that the series
  always starts at $n=0$). This is represented by
  
\centerline{\tt S=[u,h(n),[a\_1,e\_1;a\_2,e\_2;...a\_g,e\_g],w]}

  \noindent
  (when $g=0$, the empty matrix is simply set to $1$).
  Note that for convergence, we always have
  $\sum_{1\le i\le g}e_i\le0$. Thus, the above example says that
  \begin{align*}\cosh(\pi z)&=1+4z^2\sum_{n\ge0}\dfrac{(1/2-iz)_n(1/2+iz)_n}{(3/2)_n^2}\text{\quad and}\\
    \dfrac{1}{\cosh(\pi z)}&=1-4z^2\sum_{n\ge1}\dfrac{(-1/2)_n^2}{(1/2-iz)_n(1/2+iz)_n}\end{align*}
  The fact that the series begins at $n=0$ has been chosen for uniformity,
  but often leads to slightly more cumbersome formulas, for instance
  instead of $\z(2)=\sum_{n\ge1}1/n^2$ we write instead
  $\z(2)=\sum_{n\ge0}1/(n+1)^2$, hence represented by
  {\tt [0,1/(n+1)\^{}2,1,1]}.
\item {\tt v[5]} is used for testing. When empty, it means that
  the default values for the variables used in the CF can be used
  (presently $2/3$, $3/4$, $1/2$, $2$ for closures with arity less or equal to
  $4$). If it is a vector of values, these are taken instead of the default
  values. For instance for a closure of arity $2$, {\tt v[5]=[2/3,-1]}
  means that one can take the default value $2/3$ for the first variable,
  but probably cannot use $3/4$ for the second, but $-1$ instead.
  Other values of {\tt v[5]} have specific meanings which may change:
  \begin{itemize}\item {\tt v[5]=1}: Bernoulli type sequence, test expansion.
  \item {\tt v[5]=2}: Jacobi type sequence, only test expansion.
  \item {\tt v[5]=3}: Jacobi type sequence, test expansion and numerics, (2,1/2).
  \item {\tt v[5]=4}: Change z into 1/z, multiply by z, test expansion.
  \item {\tt v[5]=5}: Choose $z=2i$, multiply by $-i$, test.
  \end{itemize}
  Note that in the present version, {\tt v[5]=4} and {\tt v[5]=5} occur
  each for a single CF.
\item {\tt v[6]} is the \TeX\ label used for the CF, so as to be
  able to find it easily. It is very useful to write a trivial little
  script to find conversely the given entry of the big vector {\tt V}
  corresponding to a given label. Note that the labels are rather arbitrary
  (but of course unique) since they were given after discovery of new CFs,
  and more or less in increasing order.
\item {\tt v[7]} is a 2-component vector [A1,AD1] whose entries are strings
  describing the behavior of the CF under Ap\'ery acceleration: A1 is
  for the Ap\'ery diagonal process, and AD1 for the Ap\'ery dual. The meanings
  of the strings are as follows: ``NO'' if Ap\'ery is not applicable,
  ``BAcomp'' if it is but the Bauer--Muir method introduces complicated
  denominators, ``AP\textrightarrow label'' (``APSI\textrightarrow label'' if followed by a simplification)
  and ``APD\textrightarrow label'' if Ap\'ery leads to the corresponding labeled CF,
  (``APD\textrightarrow SELF'' it is self-dual), ``APcomp'' and ``APDcomp'' if the result is
  too complicated to be included in the encyclopedia so with no label,
  and a possible additional ``sim'' if similar to.
\item {\tt v[8]} is a (usually empty) string giving any additional information.
\end{enumerate}

\smallskip

$\bullet$ {\bf Format of Parametrized Continued Fractions}

\smallskip

Many CFs are part of families of inequivalent CFs depending on one or more
parameters. These families are included in the file in a different format,
since the initial values of $a(n)$ and $b(n)$ are not given. We have
only given such families for CFs of period 1.

The corresponding vector {\tt v} has now only five components, and again
to make the explanation clearer we give an example:

\begin{verbatim}
[[()->Pi,2*k+2,(2*n-1)*(2*n+2*k-4*u-1)]
[1,[0,-1,0,k+1]],"5.3.3",[],"u<=k/2"];
\end{verbatim}

\begin{enumerate}\item The first component {\tt v[1]} is a three-component
  vector: {\tt v[1][1]} is the closure giving the limit $S$, and {\tt v[1][2]}
  and {\tt v[1][3]} are polynomials giving $a(n)$ and $b(n)$ for $n$
  sufficiently large. This means that for suitable initial values $a(0),a(1),\dotsc$, and $b(0),b(1),\dotsc$, and for given nonnegative integers $k$ and
  $u$, the corresponding CF converges to $S$. Equivalently, the CF
  {\tt [[v[1][2]],v[1][3]]} converges to a limit of the form $(aS+b)/(cS+d)$
  with $a$, $b$, $c$, $d$ integral and $ad-bc\ne0$. For instance, choosing
  $u=0$ and $k=0,1,2$, or $u=1$ and $k=2$, the above parametrized CF says that
  the limits of the CFs
  \begin{align*}
    &2+1/(2+1/(2+9/(2+25/(2+49/(2+81/(2+\cdots))))))\;,\\
    &4-1/(4+3/(4+15/(4+35/(4+63/(4+99/(4+\cdots))))))\;,\\
    &6-3/(6+5/(6+21/(6+45/(6+77/(6+117/(6+\cdots))))))\;,\text{\quad and}\\
    &6+1/(6+1/(6+9/(6+25/(6+49/(6+81/(6+\cdots))))))\end{align*}
  are all of the form $(A\pi +B)/(C\pi + D)$ with $AD-BC\ne0$ and
  $A$, $B$, $C$, and $D$ in $\Z$.
\item {\tt v[2]} contains the information on the speed of convergence
  in the shortened format {\tt [d,[F,E,D,P]]}, since the constant {\tt C}
  and the asymptotic expansion {\tt A} vary in a parametric family so
  are not given.
\item {\tt v[3]} is the label where the parametric CF occurs (same
  as {\tt v[5]} for ordinary CFs).
\item {\tt v[4]} is used for testing, but for now is the empty vector.
\item {\tt v[5]} is a string giving some additional information.
\end{enumerate}

Note that for integer $k$, the CF $(a(n+k),b(n+k))$ is essentially the
same as the CF $(a(n),b(n))$, so is of course not considered as a new CF,
so this trivial parametrization is not included.

\smallskip

$\bullet$ {\bf Format of Definitions}

\smallskip

In some of the CFs, it is necessary to add some function definitions,
either because the corresponding functions do not exist in {\tt Pari/GP}
(for example {\tt bernfrac(k)} and {\tt eulerfrac(k)} exist, but not
{\tt tanfrac(k)} for the tangent numbers), or simply because the existing
expression is cumbersome. The command {\tt V=readvec("file.gp");} will have two
effects: first, it defines the function so that it can be used in the CFs,
and second it is kept as an element of {\tt V} as a closure (but without
its name). Thus, if you do not know the meaning of a function,
typing {\tt ?functionname} will give you an explanation or a formula for it.

\medskip

{\bf More Details}

\medskip

Here is in more detail the different types of vectors $v$ that we can
encounter, with examples.

\begin{enumerate}
\item Definitions: as explained above. These are \emph{closures}. Example:
\begin{verbatim}
  R1(k,n)=(n^k*(2*n+1)+(n-1)^k*(2*n-3))/(2*n-1);
\end{verbatim}
\item Parametrized continued fractions: as explained above. These are
  \emph{vectors} with \emph{five} components. Example:
\begin{verbatim}
[[()->2^(1/3),7*n-5+k+v+4*u,-4*(n+u)*(3*n+v-2)],
            [1,[0,4/3,0,u+(5-v)/3+2*k]],"5.2.3",[],""];
\end{verbatim}
\item Ordinary CF with $a$ or $b$ expressed as a closure: These are
  \emph{vectors} $v$ with \emph{eight} components such that $v[2][1]$
  or $v[2][2]$ is a closure. Example:
\begin{verbatim}
[k->zeta(k),[k->[0,n^k+(n-1)^k],k->[1,-n^(2*k)]],[1,[0,1,0,k-1],
k->1/(k-1),1-((k-1)/2)*x+sum(m=1,4,binomial(k+2*m-2,2*m)*
bernfrac(2*m)*x^(2*m),O(x^10))],[0,(n,k)->1/(n+1)^k,1,1],
[3],"5.7.1",["BAcomp","BAcomp"],""];
\end{verbatim}
\item Ordinary CF with the last component of $a$ and $b$ a closure. Example:
\begin{verbatim}
  [()->1,[[0,n->f(n)],[n->f(n+1)+1]],[0,0,0,0],
                              0,[],"6.1.1",["NO","NO"],""];
\end{verbatim}
\item Ordinary CF of period 2 or more: These are vectors $v$ with eight
  components such that the first (hence any) component of $v[2][1]$ or
  $v[2][2]$ is a vector. Example:
\begin{verbatim}
[()->2^(1/3),[[[0,5],[2*n,12*n+6]],[[6,-1/3],[-(n+1)*(3*n+2),
-n*(3*n+1)]]],[2,[0,(1+d)^2,0,0],()->2^(1/3)*3^(1/2)/(1+sqrt(2))^2,
1+(9*d/8)*x+(-9*d/8+81/64)*x^2+(67105*d/41472-81/32)*x^3+O(x^4)],
0,[],"5.2.10",["NO","NO"],""];
\end{verbatim}
\item Ordinary CF of period 1 with no closures. Example:
\begin{verbatim}
[()->3^(1/2),[[3/2,(2*n+1)^2],[2,-n^2*(n+2)^2]],[3,[0,(2+d)^2,0,0],
()->2*sqrt(3)/(2+sqrt(3))^3,1+d*x+(3/2-3*d/2)*x^2+(3*d-9/2)*x^3
+(45/4-27*d/4)*x^4+O(x^5)],0,[],"5.2.1",
["NO","NO"],"Infinitely contractible"];
\end{verbatim}
A second example with Ap\'ery information and series for the inverse:
\begin{verbatim}
[()->2^(1/3),[[1/2,7*n-5],[1,-12*n^2+8*n]],[1,[0,4/3,0,5/3],
()->2^(8/3)/gamma(1/3),1-55/9*x+170/3*x^2-1579490/2187*x^3
+230596751/19683*x^4+O(x^5)],[[],[0,2,[-2/3,1;1,-1],3/4]],
[],"5.2.3",["APSI->5.2.12","APD->SELF"],""];
\end{verbatim}
where the string ``APSI\textrightarrow 5.2.12'' indicates that after
simplification, Ap\'ery leads to label 5.2.12, and ``APD\textrightarrow SELF''
that it is self-dual for Ap\'ery. The series entry says that
$2^{-1/3}=2\sum_{n\ge0}((-2/3)_n/n!)(3/4)^n$,
which is simply the special case $x=3/4$ of the expansion of $2(1-x)^{2/3}$.
\end{enumerate}

\medskip

Here is a little script which outputs a number from 1 to 6 giving the type
of a vector $v$:
\begin{verbatim}
/* Find the type number of an entry in the file */

typevec(v)=
{ my(tv=type(v),lv=#v,a,b);
  if(tv=="t_CLOSURE",return(1));
  if(tv!="t_VEC"||(lv!=5&&lv!=8),error("incorrect entry"));
  if(lv==5,return(2));
  [a,b]=v[2];
  if(type(a)=="t_CLOSURE"||type(b)=="t_CLOSURE",return(3));
  if(type(a)!="t_VEC"||type(b)!="t_VEC",error("incorrect entry 2"));
  if(type(a[#a])=="t_CLOSURE"||type(b[#b])=="t_CLOSURE",return(4));
  if(type(a[1])=="t_VEC"||type(b[1])=="t_VEC",return(5));
  return(6); }
\end{verbatim}

\medskip

\subsection{Table of Contents}

\medskip

The complete GP file of the 1925 CFs and definitions together with convergence
and additional data can be obtained by obtaining the \TeX\ source and
saving the file contained between the {\tt comment} statements which
are at the end of the source file, and then reading it under {\tt GP}
with the command {\tt V=readvec("file.gp")} if {\tt file.gp} is the
name that you gave to your saved file. The source file itself is divided
into sections (included as comments so not read by {\tt GP}) which helps
to order the different CFs. The sections are the following, together
with the range of corresponding indices of the vector {\tt V}:

\bigskip

Definitions: 1--46

Constants: $2^{1/3}$: 47--86

Constants: $\pi$, $\log(2)$, and Periods of Degree 1: 87--236

Constants: $\pi^2$, $G$, and Periods of Degree 2: 237--337

Constants: $\pi^3$, $\z(3)$, and Periods of Degree 3: 338--367

Constants: $\pi^4$: 368--375

Constants: $\zeta(k)$ and $\pi^k$ for $k\ge5$: 376--379

Constants: Linear Combinations of Zeta- and $L$-Values: 380--661

Miscellaneous Constants Related to $L$-Values: 662--664

Constants Related to Powers of $e=\exp(1)$: 665--702

Constants Involving Gamma Quotients: 703--1138

Constants Coming from Bessel Functions: 1139--1158

Trivial Functions with an Arbitrary Sequence: 1159--1169

Other Trivial Functions: 1170--1182

Exponential Functions: 1183--1210

Logarithm Functions: 1211--1236

Power Functions: 1237--1262

Direct Hyperbolic and Trigonometric Functions: 1263--1332

Inverse Hyperbolic and Trigonometric Functions: 1333--1361

Gamma Function: One Variable: 1362--1416

Gamma Function: Several Variables: 1417--1432

Function $\psi(z)=\Gamma'(z)/\Gamma(z)$: 1433--1465

Function $\beta(z)=(\psi((1+z)/2)-\psi(z/2))/2$: 1466--1541

Function $\psi'(z)$: 1542--1567

Function $\beta_1(z)=-\beta'(z)=\psi'(z/2)/2-\psi'(z)$: 1568--1578

Function $\psi''(z)$: 1579--1596

Function $\beta_2(z)=-\beta'_1(z)=\psi''(z)-\psi''(z/2)/4$: 1597--1598

Function $\psi'''(z)$: 1599

Ordinary Generating Functions for Bernoulli Type Sequences: 1600--1725

Incomplete Gamma Functions: 1726--1737

Error Function and Related: 1738--1754

Jacobi Elliptic Functions: Laplace Transforms: 1755--1790

Complete Elliptic Integrals: 1791--1807

More Integrals: 1808--1830

Bessel and Related Functions: 1831--1857

General Hypergeometric Functions: ${}_0F_1$: 1858--1861

General Hypergeometric Functions: ${}_1F_1$: 1862--1869

General Hypergeometric Functions: ${}_2F_0$: 1870--1883

General Hypergeometric Functions: $U$: 1884--1885

General Hypergeometric Functions: ${}_2F_1$: 1886--1906

General Hypergeometric Functions: ${}_2F_1$ with $b=1$: 1907--1921

General Hypergeometric Functions: ${}_3F_2$ with $z=1$: 1922--1925

\medskip

Note that if you do not understand the name of a function, you can
use the {\tt GP} command {\tt ?name} to obtain an explanation or a formula.
For instance, typing {\tt ?tanfrac} returns

\begin{verbatim}
tanfrac =
  (k)->if(k==0,0,2^(k+1)*(2^(k+1)-1)*bernfrac(k+1)/(k+1))
\end{verbatim}

\medskip

{\bf Warnings.} Evidently in such a large database essentially put together
``by hand'' it is unavoidable that there are typographical
errors. But apart from this, two warnings are necessary.
\begin{enumerate}\item First, the speed of convergence has been obtained
  using heuristic methods, although I believe that the results are correct
  (in particular, for CFs for \emph{functions}, the indicated
  speed of convergence often depends on the range of the variables). In
  addition, the constant $C$ has been obtained using linear dependence
  algorithms (typically the {\tt lindep} function of {\tt Pari/GP}): it should
  be correct, but I have not tried to prove anything.
\item Second, many of the new CFs have been obtained rigorously from existing
  ones, for instance by Bauer--Muir--Ap\'ery type methods, but many others
  have been obtained ``\`a la Ramanujan'', in particular not rigorously,
  although I am quite confident of their correctness.
  Thus it is absolutely impossible for me to give the origin, and
  a fortiori the proofs, of the CFs in the database.
\end{enumerate}

Evidently I would more than welcome corrections and/or additions. Since
each CF has a label number, please refer to it when sending comments.
Ideally, this file should be put on a web site which could be regularly
updated by users, but I do not have the necessary knowledge to do this.

\medskip

\subsection{Beautifying the Continued Fractions}

\medskip

From a given CF in the database format, it is not difficult to create a
``beautified'' version which is more human-readable, and I advise the reader
to do so. In fact, the author has done so and added a large quantity of
additional information, thus creating a 676 page book, see \cite{Coh3}.
In particular, some
CFs may not be understandable without this additional information, but I
believe that the present file is sufficiently useful as is.

\smallskip

Here is an example of beautification. Entries {\tt V[1190]} and {\tt V[1191]}
are (newlines inserted for better readability):
\begin{verbatim}
[(z)->exp(Pi*z/2),[[1,-z+1,2],[2*z,4*n^2-4*n+(z^2+1)]],
[1,[0,-1,0,1],(z)->z*exp(Pi*z/2)*cosh(Pi*z/2),
1+(-1/16*z^2-1/4)*x^2+(1/128*z^4+7/64*z^2+5/16)*x^4+O(x^6)],
0,[],"6.3.8",["AP->6.3.11","APD->6.3.9"],""];

[[()->exp(Pi*z/2),4*k+2,4*n^2-4*n+(z^2+1)],
[1,[0,-1,0,2*k+1]],"6.3.8",[],""];
\end{verbatim}
A possible beautified version is:

\bigskip

{\bf 6.3.8}

\begin{verbatim}
[(z)->exp(Pi*z/2),[1,1-z,2],[2*z,z^2+(2*n-1)^2]]
\end{verbatim}
$$e^{\pi z/2}=1+\dfrac{2z}{-z+1+\dfrac{z^2+1}{2+\dfrac{z^2+9}{2+\dfrac{z^2+25}{2+\dfrac{z^2+49}{2+\dfrac{z^2+81}{2+\ddots}}}}}}$$
Convergence type $P^-$ with $P=1$ and $C=ze^{\pi z/2}\cosh(\pi z/2)$, so that
$$e^{\pi z/2}-\dfrac{p(n)}{q(n)}\sim(-1)^n\dfrac{ze^{\pi z/2}\cosh(\pi z/2)}{n}\;.$$
$$A=1+(-z^2/16-1/4)/n^2+(z^4/128+7z^2/64+5/16)/n^4+\cdots$$
Parametric family for $k\ge0$:
\begin{verbatim}
[(z)->exp(Pi*z/2),4*k+2,z^2+(2*n-1)^2]
\end{verbatim}
Convergence type $P^-$ with $P=2k+1$.

\smallskip

Ap\'ery accelerates to {\bf 6.3.11}, Ap\'ery dual is {\bf 6.3.9}.

\section{Some Remarks and Questions on the Database}

In compiling and computing this database of continued fractions, a number
of sporadic questions have arisen. Here are a few:

\begin{enumerate}
\item The CF given as an example for $3^{1/2}$ (label 5.2.1) has the
  curious property of being infinitely contractible: usually, computing
  the contraction of a CF of bidegree $(d_1,d_2)$ gives a CF of larger
  bidegree. For the present CF of bidegree $(2,4)$, the contractions
  become considerably more complicated, but can always be simplified (i.e.,
  to a CF with the same partial quotients) to bidegree $(2,4)$.
  First, I do not know how to prove this, and second, I do not know how to
  characterize such CFs.
\item The parametric families, especially those for periods of degree $2$
  and $3$ such as $\z(2)$ and $\z(3)$, have been obtained by searches.
  In many cases they are probably not complete, but in some cases there are so
  many that it would be desirable to classify them, which I have not always
  succeeded in doing. See \cite{Coh1} for detailed discussions of this.
\item In the case of logarithmic convergence, I do not know of
  any reasonable method to compute the constant $C$ such that
  $S-p(n)/q(n)\sim C/\log(n)$, except when the CF is completely explicit,
  or when it is the specialization of a CF with more variables (for instance,
  this is how I found that the constant $C$ for 5.4.27 is equal to
  $\pi^3/8$ and that of 5.5.10 is $\pi^4/28$).
\item When the convergence is precisely $S-p(n)/q(n)\sim C/n^k$ with integral
  $k$, I do not know of a systematic way to compute the asymptotic expansion
  $A$ to more than $k-1$ terms, except when the CF comes from a sum of
  explicit rational functions. In many cases, this is due to the presence of
  more complicated terms in the asymptotic expansion, typically involving
  $\log(n)$. In all other convergence cases (except of
  course the logarithmic case), including $C/n^k$ when $k\notin\Z$, there is
  an algorithmic way to compute $A$.
\item Since there are so many CFs for the functions $\be(z)$ and $\be_1(z)$,
  it is slightly surprising that there are no interesting ones for
  $\be_2(z)=-\be_1'(z)$; nonetheless, I have given a couple of (uninteresting)
  ones.
\item Although $I_{\nu}(z)$ is an even function of $z$ (up to a power of $z$),
  the CF 6.22.6 does not seem to define an even function. This implies
  that there does not seem to exist a corresponding CF for $J_{\nu}(z)$ (not
  using complex numbers), and consequently that there does not seem to
  exist analogs for $J$ of the CFs 5.12.10 and 5.12.11 for $I$.
\item There exist many CFs involving the Bessel functions $J_{\nu}(z)$,
  $I_{\nu}(z)$, and $K_{\nu}(z)$, as well as
  $H_{\nu}(z)=J_{\nu}(z)\pm i Y_{\nu}(z)$. Are there interesting ones involving
  $Y_{\nu}(z)$ alone ?
\item Using Ap\'ery, it is not difficult to find a parametric family of CFs
  for $\cosh(\pi z/3)$, see 6.6.9. On the other hand, the very similar
  CF 6.6.20 for $\sinh(\pi z/3)/(z\sqrt{3}/2)$ does not accelerate
  nicely, so I have been unable to find a corresponding parametric family.
\item The CFs 5.11.91 for $\G(1/3)^2/\G(2/3)$ and
  5.11.129 for $\G(2/3)^2/\G(1/3)$ are extremely similar and
  have the same origin. They can both be Ap\'ery accelerated giving
  complicated formulas, but the former has a very simple Ap\'ery dual
  5.11.96, while the latter has a very complicated Ap\'ery dual which
  therefore we have not given. Is there an analogue of 5.11.96 for
  $\G(2/3)^2/\G(1/3)$ ?
\item Similarly, there is a beautiful rapidly convergent CF
  5.11.140 as well as a very simple hypergeometric series for
  $\G(2/3)^2/\G(1/3)$, but apparently nothing of the sort for
  $\G(1/3)^2/\G(2/3)$.
\end{enumerate}

\vfill\eject

\end{document}